\documentclass[12pt]{article} 
\usepackage{amsmath}

\usepackage[utf8]{inputenc} 

\usepackage{geometry} 
\usepackage{graphicx} 

\usepackage{booktabs} 
\usepackage{array} 
\usepackage{paralist} 
\usepackage{verbatim} 
\usepackage{subfig} 
\usepackage{amsfonts}
\usepackage{fancyhdr} 
\usepackage{sectsty}
\allsectionsfont{\sffamily\mdseries\upshape} 

\usepackage[nottoc,notlof,notlot]{tocbibind} 
\usepackage[titles,subfigure]{tocloft} 

\begin{document}
\title{\bf {A Finite Element Scheme for an Initial Value Problem}}
\author{{\bf Vassilios K Kalpakides}\\
Department of Materials Science and Engineering,\\
University of Ioannina, Ioannina, GR-45110, Greece}
\date{} 
\maketitle
\begin{abstract}
A new Hamilton principle of convolutional type, completely compatible with the initial conditions of an IVP, has been proposed in a recent publication \cite{KALPAKIDES-et-al2019}. In the present paper the possible use of this principle for the formulation of a FE scheme adjusted to dynamical problems is investigated. To this end, a FE scheme based on a convolutional extremum principle for the harmonic oscillator -- used as an exemplary initial value problem – is developed and presented in detail.  Besides, from the local finite element analysis a recurrent (one-step) algorithm arises which provides an approximate solution to the IVP, as well. The succeeded schemes are computationally tested for both free and forced vibration problems.
\end{abstract}
\section{Introduction}
It is widely accepted that the initial value problems do not admit appropriate variational principles like the corresponding boundary value ones.
This is due to the fact that the differential operator of dynamical problems, on account of the initial conditions, are not symmetric. If one considers solely the differential operator, independently of the initial conditions, it might  be symmetric in some cases. Then it will be referred to as {\em formally symmetric} \cite{TONTI1973} and the formulation of a variational principle of Hamilton's type is possible.  However, the initial conditions can not emerge from the variational principle and must be artificially  imposed on the differential equation. \\
\\
As Tonti remarked in \cite{TONTI1973}, the notion of symmetry is defined  with respect to a bilinear form, often hidden in the background. Apparently, for differential operators, the symmetry is defined with respect to the  $L^2$ inner product. Thus more precisely, one can state that  {\em the operators of the initial value problems are not symmetric with respect to the symmetry rising from the  $L^2$ inner product}. In other words,  symmetry is not an inherent property of an operator; it should be regarded with respect to  any specific bilinear form. An operator which is not symmetric with respect to the inner product, might be symmetric with respect to another bilinear form. \\
\\
 From this viewpoint, Gurtn's attempt \cite{GURTIN1963,GURTIN1964,GURTIN1964b} to formulate variational pronciples for viscoelasticity and  elastodynamics using the convolution  was quite ineresting.  The  same  approach, a littlle improved, was followed Reddy \cite{REDDY1976, REDDY1976b} and later by a great number of authors  \cite{APOSTOL-et-al2013, APOSTOL-et-al2013b,DARGUSH2012,DARGUSH-et-al2012,DARGUSH-et-al2015,DARGUSH-et-al2016,KIM2014,PENG-et-al1996,RAFAL1968} establishing so what one can call {\em convolutional variatonal principles} \cite{AMIRI-et-al2029}. Actually, in these works the new variational principles have been obtained by the use of the convolution bilinear form instead of the standard one induced by $L^2$ inner product.\\
\\
Once a convolutional  variational principle for an initial value problem has been obtained, one may wonder whether or not it is possible to develop a finite element scheme analogous to the standard one for  boundary value problems. Dargush \cite{DARGUSH2012} and his co-workers \cite{KIM2014,DARGUSH-et-al2015,DARGUSH-et-al2016,DARRALL-et-al2018} have presented quite interesting work  in this direction. They  exploited their mixed convolutional  variational principles  to conclude in FE approximation schemes for initial--boundary value problems. However,   as concerns time evolution, these schemes are  one--step marching algorithms and, to our best knowledge, the question on the possibility of a true FE scheme for the unified space--time continuum  remains unanswered.\\
\\
The present work is a contribution towards answering this question. Here the posibility for the formulation  of a true FE scheme for the time continuum is proved, in other words, a formulation of a FE scheme for an initial value problem is presented.  Following up an earlier work \cite{KALPAKIDES-et-al2019}, where  a generalized Hamilton's  principle in a convolutional form had been proposed, we present a fully  developed finite elment formulation which provides an aproximate solution to the variational equation. A second approximate scheme is produced as a by-product of the FE formulation process. The local  algebraic system that is formulated for the needs of the FE analysis provides a recurrent formula, like that one of \cite{DARGUSH2012}, upon which a one--step algorithm is fully  developed, as well.
It is worth remarking  that the adopted convolutional symmetry is "transferred" to the finite dimensional version of the problem, where it  appears  as  symmetry with respect to the second diagonal of the global stiffness matrix.\\
\\
In Section 2, some basic definitions and propositions  are recalled and wherever necessary are proved. A new definition for  the convolution at a subdomain (element) as well as a corresponding variational principle at element level are provided. Also, the partition of the initial variational problem into a finite number of variational problems, one for each element, is presented. In Section 3, the FE scheme is analytically developed. First, a typical FE analysis is presented; starting from the element level, the local systems are formulated  and appropriately assembled to the global algebraic system. Next, in Section 4, the inverse procedure is examined, by inserting a linear combination of the global shape functions into the variational equation, concluding in the same algebraic system. In Section 5, an alternative, one step approximation scheme stemming from the local algebraic system is presented and its theoretical critical value for numerical stability is provided. Finally, in Section 6, some computational examples are presented  supporting the  theoretical considerations and the efficiency  of the two numerical schemes to approximate the analytical solution. The paper ends up with Section 7, where the main conclusions are stated.
\section{Prerequisites}
\subsection{Definitions and propositions}
In this paper, the simple problem of an harmonic oscillator is used as an example  of an initial value problem. Thus, both theory and numerical examples concen the following problem 
\\
\\
{\em {\bf Problem 1}. Find a $u$ in  $C^2[0,t]$ that satisfies the differential equation 
\begin{equation}
mu''(s)+ku(s)=f(s),\   s\in (0,t) , \ t>0
\end{equation}
and the initial conditions}
\begin{equation}
u(0)=u_0, \quad
 u'(0)=v_0,
\end{equation}
\\
where  $f\in C[0,\infty)$ is a periodic  function, $m$ and $k$ are positive constants, $u_0$ and $v_0$ are the initial  data.
\\
Because the notion of convolution  is of central importance in our analysis, we recall its definition\\
\\
 {\bf Definition 1} {\em If  $f, \ g$ are $L^2$ functions over $(0,t)$, then their convolution is defined as}
\begin{equation}
[g,h](t)=\int_0^tg(s)h(t-s)ds,\quad t>0.
\end{equation}
In what follows, if there is no danger of confusion, the argument in the above definition will be omitted.\\
\\
Discretizing the interval $[0,t]$ in finite subintervals, a convolution over any finite interval of the form $[t_1,t_2],\ t_1>0$ will be necessary.  For that purpose the following definition is proposed\\
\\
{\bf Definition 2} {\em Let $[t_1,t_2]$ any finite  interval of  $\mathbb{R}$ with $t_1\geq 0$. The convolution of the  functions  $g, h \in L^2(t_1,t_2)$ is defined as}
\begin{eqnarray}
[g,h]_{t_1}^{t_2}&=&\int_0^{\tau}g(t_1+s)h(t_2-s)ds\nonumber\\
&=&\int_{t_1}^{t_2}g(s)h(t_1+t_2-s)ds, \quad \tau=t_2-t_1.
\end{eqnarray}
\\
\\
It is easy to confirm that the above definition keeps all desired properties of the standard definition and it falls to it  when $t_1=0$.
\\
\\
It has been proved \cite{KALPAKIDES-et-al2019} that the solution of the initial value problem (1)--(2) is related to  the problem of finding a statioanary value for the functional
\begin{equation}
I[u]=\frac{m}{2}\left[ u',u'\right]+\frac{k}{2}\left[ u,u\right]-[f,u]-[\tilde{f},u],\quad u\in { \cal D},
\end{equation}
where $\tilde{f}(s)=f_0\delta (s),\ f_0=mv_0$,  $\delta$  is the Dirac function  and 
\begin{equation}
 {\cal D}=\{\phi \in C^2[0,t]\  \text{and} \  \phi(0)=u_0\}.
\end{equation}
More specifically, it holds\\
\\
{\bf Proposition 1} {\em If the functional $I$ attains a statioanry value at some $u\in {\cal D}$, then that  $u$ will be a  solution to the initial value problem (1)--(2), too}.\\
\\
An analogous propostion holds for an interval of the form $[t_1,t_2]$.  Consider the \\
\\
{\bf Problem 2} {\em Find a function  $u$ in $C^2[t_1,t_2]$ which satisfies the  following differential equation and initial conditions}
\begin{equation}
mu''(s)+ku(s)=f(s),\   s\in (t_1,t_2) , \ t_1>0,
\end{equation}
\begin{equation}
u(t_1)=u_1,   u'(t_1)=v_1.
\end{equation}
\\
This problem is related to the functional 
\begin{equation}
I^t[u]=\frac{m}{2}\left[ u',u'\right]_{t_1}^{t_{2}}+\frac{k}{2}\left[ u,u\right]_{t_1}^{t_2}-[f,u]_{t_1}^{t_2}-[\tilde{f^1},u]_{t_1}^{t_2},\quad
 u\in{\cal D}^t,
\end{equation}
where
\begin{equation}
{ \cal D}^t=\{\phi \in C^2[t_1,t_2]\  \text{and} \  \phi(t_1)=u_1\}
\end{equation}
and 
 $$
\tilde{f^1}(s)=f^1\delta(s-t_1), \quad f^1=mv_1.
$$
The following statement determines the relation between Problem 2 and the functional $I^t$\\
\\
{\bf Proposition 2} {\em If the functional $I^t$ attains a stationary value at some $u\in {\cal D}^t$, then that  $u$ will be a  solution to the  problem (7)--(8), too}.\\
\\
{\bf Proof} The functional $I^t$ is transformed by translating all functions within its argument  to the left by an interval $t_1$
$$
h(s)=h^*(s-t_1).
$$
For instance, the  second term of the functional
\begin{equation*}
\left[u,u\right]_{t_1}^{t_2}=\int_{0}^{t_2-t_1}u(t_1+s)u(t_2-s)ds
\end{equation*}
becomes
$$
\left[u,u\right]_{t_1}^{t_2}=\int_0^{t_2-t_1}u^*(s)u^*(t_2-t_1-s)ds=\int_0^{\tau}u^*(s)u^*(\tau-s)ds=\left[u^*,u^*\right](\tau),
$$
where $\tau =t_2-t_1$.\\
Also, the last term becomes
\begin{eqnarray*}
[\tilde{f^1},u]_{t_1}^{t_2}=\int_0^{t_2-t_1}\tilde{f^1}(t_1+s)u(t_2-s) ds=\int_0^{\tau}\tilde{f^1}^*(s)u^*(\tau-s) ds=\left[\tilde{f^1}^*,u^*\right]\nonumber \\
=\int_0^{\tau}f^1\delta(s)u^*(\tau-s) ds=f^1u^*(\tau)
\end{eqnarray*}
Similarly, accounting for the rest terms of $I^t$, one concludes
\begin{equation}
I^t[u]=I^*[u^*]=\frac{m}{2}\left[{ u^*}',{u^*}'\right]+\frac{k}{2}\left[ u^*,u^*\right]-\left[f^*,u^*\right]-\left[{\tilde{f^1}}^*,u^*\right],\quad
 u\in{{\cal D}^t}^*,
\end{equation}
where
$$
{{ \cal D}^t}^*=\{\phi \in C^2[0,t]\  \text{and} \  \phi(0)=u_1\},
$$
Notice also that 
$$
\tilde{f^1}^*=f^1\delta(s),\ f^1=mv_0
$$
Thus, if $I^*$ attains a stationanry value at $u^*$, according to Proposition 1, it will satisfy eq. (1) and the initial conditions (2), i.e.,
\begin{equation*}
m{u^*}''(z)+ku^*(z)=f^*(z),\   z\in (0,\tau) , \ \tau>0,
\end{equation*}
\begin{equation*}
u^*(0)=u_0,  { u^*}'(0)=v_0,
\end{equation*}
or, equivalently on account of eq. (12)
\begin{equation*}
mu''(z+t_1)+ku(z+t_1)=f(z+t_1),\   z\in (0,\tau) , \ \tau>0,
\end{equation*}
\begin{equation*}
u(t_1)=u_0,  u'(t_1)=v_0.
\end{equation*}
 Introducing the substitution
$$s=z+t_1$$
the above equations can be written
\begin{equation*}
mu''(s)+ku(s)=f(s),\   s\in (t_1,t_2) ,
\end{equation*}
\begin{equation*}
u(t_1)=u_0,  u'(t_1)=v_0
\end{equation*}
and the proof has been completed.\\
\\
Finally, based on the above  proposiitons, one can prove  that the following statement holds\\
\\
{\bf Proposition 3} {\em Consider  the functional 
\begin{equation}
I^{\tau}[u]=\frac{m}{2}\left[ u',u'\right]_{t_1}^{t_{2}}+\frac{k}{2}\left[ u,u\right]_{t_1}^{t_2}-[f,u]_{t_1}^{t_2}-[\tilde{f^1}-\tilde{f^2},u]_{t_1}^{t_2},\quad
 u\in {\cal D}^{\tau}=C^2[t_1,t_2],
\end{equation}
where  $\tilde{f^1}=f^1\delta(s-t_1)$,  $\tilde{f^2}=f^2\delta(s-t_2)$  and  $f^1$ , $f^2$  are given constants.
If the functional $I^{\tau}$ attains a stationary value at some  $u\in{\cal D}^{\tau}$, then that $u$ will be the solution of the following problem}
\begin{equation}
mu''(s)+ku(s)=f(s),\   s\in (t_1,t_2) , 
\end{equation}
\begin{equation}
u'(t_1)=f^1/m,\quad u'(t_2)=f^2/m.
\end{equation}
\\
{\bf Remark 1} The problem (13)--(14) is not well--posed.  Notice that the data in eq. (14) concern both ends of the time interval $[t_1,t_2],$ a fact that is not compatible with an initial value problem. In addittion, the values of the velocity at $t_1$ and $t_2$ are generally unknown. Later on, it will become apparent that these values are not necessary for the numerical solution of the problem.\\
\\
{\bf Remark 2} As concerns the smoothness of the solutions in the above propositions, notice that it was required more smoothness so as to establish equivalence between the variational and the classical problem. If one is interested only in the variational problem, then requiring a solution in $C^1[0,t]$ is enough. Moreover, if  the weak form of the problem is of interest,  $W^1_2(0,t)$  should be the  appropriate space for a   generalized solution. 
\\
\subsection{Partition of the Problem}
To deal with an approximate solution of the Problem 1, we apply a partition of the domain $[0,t]$ in $n$ subintervals.
$$
0=s_1<s_2<\dots<s_n<s_{n+1}=t
$$
Though  elements of equal length are usually sufficient  for a time discretization, one can take randomly the node locations provided that the number of elements is even  ($n=2k,\ k\in \mathbb{Z}$) and the distribution of nodes in $[0,t]$ is  symmetric with respect to the midle point of the interval. That means,  the $k$ node  is necessarily located at the point $t/2$, the $i-$node is symmetric to $(n+2-i)-$node ($i=1,\dots,n+1$) and the $e-$element is symmetric and of equal length to $(n+1-e)-$element ($e=1,\dots,n$).\\
\\
Let  $e$ be any element but the first one ( $e>1$) with end points $s_e$ and $s_{e+1}$, respectively. Since there are no data for this element, to formulate a variational problem, we take the functional
$$
I^e[u]=\frac{m}{2}\left[ u',u'\right]_{s_e}^{s_{e+1}}+\frac{k}{2}\left[ u,u\right]_{s_e}^{s_{e+1}}-[f,u]_{s_e}^{s_{e+1}}-\left([\tilde{f^e_1},u]_{s_e}^{s_{e+1}}-[\tilde{f^e_2},u]_{s_e}^{s_{e+1}}\right), \quad u\in{\cal D}^e,
$$
where
$$
 {\cal D}^e= C^2[s_e,s_{e+1}]
$$
and
 $$
\tilde{f^e_1}=f^e_1\delta(s-s_e),\quad  \tilde{f^e_2}=f^e_2\delta(s-s_{e+1})
$$
with $f^e_1$ and $f^e_2$  are unknown  constants.\\
\\
If the functional $I^e$ takes  a stationary value at  $u^e\in{\cal D}^e$, then that $u^e$ will be a solution to the problem
\begin{eqnarray}
mu''(s)+ku(s)=f(s),\   s\in (s_e,s_{e+1}) , \nonumber\\
\\
u'(s_e)=f^e_1/m,   u'(s_{e+1})=f^e_2/m.\nonumber
\end{eqnarray}
Thus, the variational problems
$$
I^e[u],\  u\in{\cal D}^e \rightarrow \text{stationary}
$$
 are linked with the problems (15), for all e with  $2\le e\le n$.
\\
\\
The first element, $[s_1,s_2]$, satisfies the initial conditions of the Problem 1, that is
$$
u(s_1)=u_0,\ u'(s_1)=v_0,
$$
which according to Proposition 2, leads to the variational problem
$$
I^1[u],\  u\in{\cal D}^1 \rightarrow \text{stationary},
$$ 
where
\begin{equation*}
I^1[u]=\frac{m}{2}\left[ u',u'\right]_{s_1}^{s_2}+\frac{k}{2}\left[ u,u\right]_{s_1}^{s_2}-[f,u]_{s_1}^{s_2}-[\tilde{f^1_1},u]_{s_1}^{s_2},
\end{equation*}
with 
$$
\tilde{f^1_1}=mv_1\delta(s-s_1)
$$ 
and  
$$
 {\cal D}^1=\{\phi \in C^2[s_1,s_2]\  \text{and} \  \phi(s_1)=u_1\}.
$$
\\
\\
In conclusion,  the variational  form of  Problem 1, 
\begin{equation}
I[u],\  u\in{\cal D} \rightarrow \text{stationary},
\end{equation}
can be  divided in $n$  variational problems, one for each  element:
\begin{equation}
\begin{array}{c}
I^1[u],\  u\in{\cal D}^1 \rightarrow \text{stationary},\\
 I^2[u],\  u\in{\cal D}^2 \rightarrow \text{stationary},\\
\vdots\\
I^n[u],\  u\in{\cal D}^n \rightarrow \text{stationary}.
\end{array}
\end{equation}
Thus, the variational problem (16) is equivalent to the collection of the variational problems (17).
\section{A Finite Element Scheme}
Looking for an approximate solution of the Problem 1, we relax the requirements for the smoothness of the solution to the variational problem, allowing for piecewise linear approximations. Thus, at element level, approximations in the form of a polynomial of first order are sought. Also, for the moment, we are not concerned with the initial conditions of Problem 1. These will be imposed at the last stage of the analysis, when the global algebraic system will be established.\\
Summing up, for approximating the solution we consider the collection of variational problems:
\begin{eqnarray}
I^e[u] \rightarrow \text{stationary} \ \text{in}\  X^e,\quad e=1,\dots,n,
\end{eqnarray} 
where $X^e$  is the space of first order polynomials over $[s_e,s_{e+1}].$
\subsection{Local Analysis}
Let any element $e$ of the partition limited by the nodes $e$  and  $e+1$. Looking for a stationary point in $X^e$ for the functional $I^e$, one may consider  an approximation of the form
\begin{equation}
u^e(s)=u^e_1N^e_1(s)+u^e_2N^e_2(s),
\end{equation}
where 
\begin{equation}
N^e_1(s)=\frac{s_{e+1}-s}{s_{e+1}-s_e}, \quad N^e_2(s)=\frac{s-s_e}{s_{e+1}-s_e},\quad s\in [s_e,s_{e+1}]
\end{equation}
and $u_1^e$, $U_2^e$ are constants to be determined.\\

Inserting the solution $u^e$ into the functional $I^e$, one obtains
\begin{eqnarray}
I^e[u^e]=\frac{m}{2}\left[{u^e}',{u^e}'\right]^{s_{e+1}}_{s_e}+\frac{k}{2}\left[u^e,u^e\right]^{s_{e+1}}_{s_e}-\left[f,u^e\right]^{s_{e+1}}_{s_e}
-\left[\tilde{f^e_2},u^e\right]^{s_{e+1}}_{s_e}\nonumber\\
+\left[\tilde{f^e_1},u^e\right]^{s_{e+1}}_{s_e}
\end{eqnarray}
After a long calculation (see Appendix A), eq. (21)  may be written in the form of  a quadratic function of $u_i^e$:
\begin{eqnarray*}
\hat{I}^e(u^e_i)&=&\frac{1}{2}\left(M^e_{ij}u^e_iu^e_j+K^e_{ij}u^e_iu^e_j\right)-F^e_iu^e_i-f^e_2u^e_1+f^e_1u^e_2, \nonumber\\
&=&\frac{1}{2}{\cal K}^e_{ij}u^e_iu^e_j-F^e_iu^e_i-f^e_2u^e_1+f^e_1u^e_2,\\
&&\hspace{5cm}\quad e=1,\dots,n,\quad i,j=1,2 \nonumber
\end{eqnarray*}
or, in matrix denotation
\begin{equation}
\hat{I}^e({\bf u}^e)=\frac{1}{2}{{\bf u}^e}^T{\cal K}^e{\bf u}^e-{{\cal F}^e}^T{\bf u}^e,
\end{equation}
where 
\begin{equation}
{\cal K}^e={\bf M}^e+{\bf K}^e,\quad {\cal F}^e={\bf F}^e+{\bf f}^e
\end{equation}
are the {\em local stiffness matrix} and the {\em local force vector}, respectively.\\
The matrices ${\bf M}^e$,  ${\bf K}^e$ and the vectors ${\bf F}^e$, ${\bf f}^e$ are  given by the relations\\
\begin{equation}
M^e_{ij}=m\left[{N^e_i}',{N^e_j}'\right]^{s_{e+1}}_{s_e},\quad
K^e_{ij}=k\left[{N^e_i},{N^e_j}\right]^{s_{e+1}}_{s_e},
\end{equation}
\begin{equation}
F^e_i=\int_0^{\tau_e}f(s_e+s)N^e_i(s_{e+1}-s)ds,
 \quad
{\bf f}^e=\begin{bmatrix} f^e_2\\-f^e_1 \end{bmatrix}=\begin{bmatrix} mu'(s_{e+1})\\-mu'(s_e) \end{bmatrix}.
\end{equation}
for $i,j=1,2,\ e=1,\dots,n$.\\
\\
A necessary condition for the function $\hat{I}^e$ to take a stationary value at ${\bf u}^e$ is the satisfaction of the following relations
$$
\frac{\partial\hat{I}^e}{\partial u^e_1}=0,\quad  \frac{\partial \hat{I}^e}{\partial u^e_2}=0,
$$
or, with the aid of eq. (22)
\begin{equation} {\cal K}^e{\bf u}^e={\cal  F}^e,
\end{equation}
for $e=1,\dots,n.$\\
\\
{\bf Reamark 3} Notice that, due to the symmetry of the node locations that has been adopted, the local stiffness matrices enjoy the equalities
\begin{equation}
{\cal K }^e={\cal K}^{n+1-e}, \quad e=1,\dots,n.
\end{equation}
{\bf Remark 4} Also with the aid of Definition 2 and eq. (20), one can easily verify the following useful relations
\begin{eqnarray}
\left[f,N_1^e\right]_{s_e}^{s_{e+1}}=\int_0^{\tau_e}f(s_e+s)N_1^e(s_{e+1}-s)ds=\int_{s_e}^{s_{e+1}}f(s)N_2^e(s)ds=\left(f,N_2^e\right),\nonumber \\
\\
\left[f,N_2^e\right]_{s_e}^{s_{e+1}}=\int_0^{\tau_e}f(s_e+s)N_2^e(s_{e+1}-s)ds=\int_{s_e}^{s_{e+1}}f(s)N_1^e(s)ds=\left(f,N_1^e\right),\nonumber
\end{eqnarray}
where the parenthesis with comma $(\ , \ )$ denotes the standard inner product of $L^2$.
\subsection{The Assembly}
Now, one may unite the solutions at the different elements, accounting for the fact that the displacements at the common nodes of the neighbouring elements should coincide, ensuring the continuity of the approximate solution.\\
\\
Introducing global enumeration for the dispacements 
\begin{equation}
U_1=u^1_1,\ U_{i+1}=u^i_2=u^{i+1}_1,\ U_{n+1}=u^n_2, \quad i=1,\dots.n-1,
\end{equation}
the local systems provided by eqs. (26) are written in extensive form
\begin{eqnarray*}
{\cal K}^1_{11}U_1+{\cal K}^1_{12}U_2=F^1_1-mu'(s_2),\nonumber\\
{\cal K}^1_{21}U_1+{\cal K}^1_{22}U_2=F^1_2+mu'(s_1),\nonumber\\
\nonumber\\
{\cal K}^2_{11}U_2+{\cal K}^2_{12}U_3=F^2_1-mu'(s_3),\nonumber\\
{\cal K}^2_{21}U_2+{\cal K}^2_{22}U_3=F^2_2+mu'(s_2),\nonumber\\
\dots\dots\dots\dots\dots\dots\dots\dots\dots\dots.\nonumber\\
\dots\dots\dots\dots\dots\dots\dots\dots\dots\dots.\nonumber\\
{\cal K}^{n-1}_{11}U_{n-1}+{\cal K}^{n-1}_{12}U_n=F^{n-1}_1-mu'(s_n),\\
{\cal K}^{n-1}_{21}U_{n-1}+{\cal K}^{n-1}_{22}U_n=F^{n-1}_2+mu'(s_{n-1}),\\
\nonumber\\
{\cal K}^n_{11}U_n+{\cal K}^n_{12}U_{n+1}=F^n_1-mu'(s_{n+1}),\nonumber\\
{\cal K}^n_{21}U_n+{\cal K}^n_{22}U_{n+1}=F^n_2+mu'(s_n),\nonumber
\end{eqnarray*}
Adding those equations that have  the last term in common, one obtains
\begin{eqnarray*}
{\cal K}^1_{21}U_1+{\cal K}^1_{22}U_2=F^1_2+mu'(s_1),\\
{\cal K}^1_{11}U_1+\left({\cal K}^1_{12}+{\cal K}^2_{21}\right)U_2+{\cal K}^2_{22}U_3=F^1_1+F^2_2,\\
\dots\dots\dots\dots\dots\dots\dots\dots\dots\dots.\nonumber\\
\dots\dots\dots\dots\dots\dots\dots\dots\dots\dots.\nonumber\\
{\cal K}^{n-1}_{11}U_{n-1}+\left( {\cal K}^{n-1}_{12}+{\cal K}^n_{22} \right) U_n+{\cal K}^n_{22}U_{n+1}=F^{n-1}_1+F^n_2,\\
{\cal K}^n_{11}U_n+{\cal K}^n_{12}U_{n+1}=F^n_1-mu'(s_{n+1}),
\end{eqnarray*}
Thus, a system of $n+1$ equations has arisen. The next step is to inverse  (turn upside down) the sequence of the equations
\begin{eqnarray*}
{\cal K}^n_{11}U_n+{\cal K}^n_{12}U_{n+1}=F^n_1-mu'(s_{n+1}),\\
{\cal K}^{n-1}_{11}U_{n-1}+\left( {\cal K}^{n-1}_{12}+{\cal K}^n_{22} \right) U_n+{\cal K}^n_{22}U_{n+1}=F^{n-1}_1+F^n_2,\\
\dots\dots\dots\dots\dots\dots\dots\dots\dots\dots.\nonumber\\
\dots\dots\dots\dots\dots\dots\dots\dots\dots\dots.\nonumber\\
{\cal K}^1_{11}U_1+\left({\cal K}^1_{12}+{\cal K}^2_{21}\right)U_2+{\cal K}^2_{22}U_3=F^1_1+F^2_2,\\
{\cal K}^1_{21}U_1+{\cal K}^1_{22}U_2=F^1_2+mu'(s_1),
\end{eqnarray*}
Recalling the element symmetry given by eq. (27), it is apparent that it holds  
$$
{\cal K}^n={\cal K}^1,\quad  {\cal K}^{n-1}={\cal K}^2,\quad  \dots
$$ 
and the system of algebraic equations  takes its final form
\begin{eqnarray}
{\cal K}^1_{11}U_n+{\cal K}^1_{12}U_{n+1}=F^n_1-mu'(s_{n+1}),\\
{\cal K}^2_{11}U_{n-1}+\left( {\cal K}^2_{12}+{\cal K}^1_{22} \right) U_n+{\cal K}^1_{22}U_{n+1}=F^{n-1}_1+F^n_2,\\
\dots\dots\dots\dots\dots\dots\dots\dots\dots\dots.\nonumber\\
\dots\dots\dots\dots\dots\dots\dots\dots\dots\dots.\nonumber\\
{\cal K}^n_{11}U_1+\left({\cal K}^n_{12}+{\cal K}^{n-1}_{21}\right)U_2+{\cal K}^{n-1}_{22}U_3=F^1_1+F^2_2,\\
{\cal K}^n_{21}U_1+{\cal K}^n_{22}U_2=F^1_2+mu'(s_1).
\end{eqnarray}
The system of eqs. (30)--(33) is  written in matrix form
\begin{equation*}
\begin{bmatrix}
0 & \dots &\dots & 0 & {\cal K}^1_{11}& {\cal K}^1_{12}\\
0 & \dots & 0 & {\cal K}^2_{11}&  {\cal K}^2_{12}+  {\cal K}^1_{21}& {\cal K}^1_{22}\\
\vdots &\vdots   & \vdots &\vdots & \vdots & \vdots\\
 {\cal K}^n_{11} & {\cal K}^n_{12} + {\cal K}^{n-1}_{21} & {\cal K}^{n-1}_{22} &0&\dots&0\\
 {\cal K}^n_{21} & {\cal K}^n_{22} & 0 &\dots&\dots&0
\end{bmatrix}
\begin{bmatrix}
U_1\\U_2\\ \vdots\\U_n\\ U_{n+1}
\end{bmatrix}
=\begin{bmatrix}
{\cal F}_1\\{\cal F}_2\\ \vdots\\{\cal F}_n\\ {\cal F}_{n+1}
\end{bmatrix}
\end{equation*}
or, 
\begin{equation}
{\cal K}{\bf U}={\cal F}
\end{equation}
where
\begin{eqnarray}
{\cal K}&=&
\begin{bmatrix}
{\cal K}_{11} &{\cal K}_{12} & \dots& {\cal K}_{1n} & {\cal K}_{1(n+1)}\\
{\cal K}_{21} &{\cal K}_{22} & \dots& {\cal K}_{2n} & {\cal K}_{2(n+1)}\\
\vdots & \vdots&\ddots &\vdots&\vdots \\
{\cal K}_{n1} &{\cal K}_{n2} & \dots& {\cal K}_{nn} & {\cal K}_{n(n+1)}\\
 {\cal K}_{(n+1)1}  & {\cal K}_{(n+1)2}&\dots & {\cal K}_{(n+1)n}& {\cal K}_{(n+1)(n+1)}
\end{bmatrix} \nonumber\\
\nonumber \\
&=&
\begin{bmatrix}
0 & \dots &\dots & 0 & {\cal K}^1_{11}& {\cal K}^1_{12}\\
0 & \dots & 0 & {\cal K}^2_{11}&  {\cal K}^2_{12}+  {\cal K}^1_{21}& {\cal K}^1_{22}\\
\vdots &\vdots   & \vdots &\vdots & \vdots & \vdots\\
 {\cal K}^n_{11} & {\cal K}^n_{12} + {\cal K}^{n-1}_{21} & {\cal K}^{n-1}_{22} &0&\dots&0\\
 {\cal K}^n_{21} & {\cal K}^n_{22} & 0 &\dots&\dots&0
\end{bmatrix}
\end{eqnarray}
and
\begin{equation}
{\cal F}=
\begin{bmatrix}
{\cal F}_1\\{\cal F}_2\\ \vdots\\{\cal F}_n\\ {\cal F}_{n+1}
\end{bmatrix}
=\begin{bmatrix}
F^n_1\\F^n_2+F^{n-1}_1\\ \vdots\\F^1_1+F^2_2\\F^1_2
\end{bmatrix}
+\begin{bmatrix}
-m u'(s_{n+1})\\ 0\\ \vdots\\ 0\\ mu'(s_1)
\end{bmatrix}
\end{equation}
are the {\em  global stiffness matrix} and the {\em global force vector}, respectively. \\
\\
{\bf Remark 5} It is noted that the stiffness matrix, as it appears in eq. (35), does not take the typical form of the standard finite element method. Here, the non-vanishing entries lie along the second diagonal of the matrix. Actually, the notion of the matrix symmetry has been changed and now is regarded with respect to the second diagonal. The standard matrix symmetry comes out from the common scalar product of the finite dimensionsl Euclidean spaces. Here an alternative symmetry arises behind  which is hidden the  bilinear form 
$$
B({\bf x},{\bf y})=x_1y_{n+1}+x_2y_{n}+\dots+x_ny_2+x_{n+1}y_1,\quad {\bf x},{\bf y}\in \mathcal{E}^{n+1}
$$
which is compatible with the convolution upon which the present analysis is based.\\
\\
{\bf Remark 6} The loading vector $\cal F$  is noticeble, too.  The  initial condition (2b) which applies at the first node ($s_1=0$) is pushed away from  it and contributes to the last node ($s_{n+1}=t$). Certainly,  the  initial condition (2a) contributes to the first node through the set of admissible functions $\cal D$ (see eq. (6)).  This  redistribution of the two initial conditions  to the initial and final time   is essentialy the  mechanism that restores the inconsistency of the standard variational principles with the initial conditions. 
\subsection{Imposition of the first  Initial Condition}
At this stage one may impose the first initial condition of the problem. The condition (2a) provides the displacement at the first node, i.e.,  $U_1=u(0)=u_0$. One may deal with this condition as one does with the Dirichlet conditions in a boundary value problem.  Thus one may ignore the first equation of the system (34), or equivalently may omit the first row and the first column of the golabal stiffness matrix and modify appropriately the global force vector.\\
\\
As concerns the second condition (2b), it has been included in the functional $I^1$, (see eqs.  (21) and (25b))  resulting in the appearance of the  initial  momentum $mu'(s_1)$  in the last row of the global force vector. In conclusion, one may establish an $n \times n$ stiffness matrix and a fully known $n-$dimensional  force vector to determine the $n$ unknown parameters $U_2,\dots,U_{n+1}$ .
\section{Global Analysis}
In the previous section we started, as usual, from the local stiffness matrix and the corresponding linear system at each element and assembling we conclude in the global system. Because the presented approximation scheme is not the standard FEM, we  will reproduce this result starting from the global form of the sought solution.\\
\\
The global finite element solution over the entire $[0,t]$ will be taken  as a linear combination of the global shape functions which  are composed from the local ones as
\begin{equation}
N_i(s)=\left\{\begin{array}{cc}
N_2^{i-1}(s),&s\in[s_{i-1},s_i]\\
N_1^i(s), &s\in[s_{i},s_{i+1}]\\
0,&\text{elsewhere},
\end{array}
\right.
\quad   i=1,\dots,n+1.
\end{equation}
The functions $N_i$ enjoy the relation 
$$
N_i(s_j)=\delta_{ij},
$$
i.e, are linear independent  and span the finite dimensional space of the piecewise linear functions over $[0,t]$:
$$
X_{n+1}:=\text{span}\{N_1,\dots,N_{n+1}\}\subset W^1_2(0,t).
$$
Thus, according to Proposition 1, to find  an approximate solution of the Problem 1, in the  form of a linear combination of the global shape functions is equivalent to the\\
\\
{\bf Approximation problem} {\em Find  a solution to the variational problem
\begin{equation}
I[u]\rightarrow \text{stationary},
\end{equation}
in $X_{n+1}$.}\\
\\
Consider now such an approximate solution in the form of
\begin{equation}
u^\text{fe}(s)=\sum_{i=2}^{n+1}U_iN_i(s),\quad s\in [0,t]
\end{equation}
where $U_i\in \mathcal{R},\ i=2,\dots n+1$ are the unknown displacements that has to be computed. Notice that the first shape function $N_1$ has been omitted, because the displacement $U_1$ is known by the first initial condition  (see eq. (2a)).\\
\\
It will be shown that the unknown parameters $U_i$ can be determined by solving the linear system (34).
\\
\\
Indeed, inserting eq.  (39)   into the functional  $I$  (given by eq. (5)), one obtains
\begin{equation}
I[u^{fe}]=\hat{I}(U_i)=\frac{1}{2}M_{ij}U_iU_j+\frac{1}{2}K_{ij}U_iU_j-F_iU_i-\tilde{F}_iU_i,
\end{equation}
or
\begin{eqnarray*}
\hat{I}({\bf U})&=&\frac{1}{2}{{\bf U}}^T{\bf M}{\bf U}+\frac{1}{2}{{\bf U}}^T{\bf K}{\bf U}-{{\bf F}}^T{\bf U}-{\tilde{\bf F}}^T{\bf U}\\
&=&\frac{1}{2}{{\bf U}}^T{\cal K}{\bf U}-{{\cal F}}^T{\bf U}
\end{eqnarray*}
where
\begin{equation}
{\bf {\cal K}=M+K,\quad {\cal F}=F+\tilde{F}}.
\end{equation}
and
\begin{eqnarray}
M_{ij}=m[N'_i,N'_j],\  K_{ij}=k[N_i,N_j],\ F_i=[f,N_i],\ \tilde{F}_i=[\tilde{f},N_i],\\
  i,j=2,\dots,n+1.\nonumber
\end{eqnarray}
A necessary condtion for $\hat{I}$ to attain a stationary value at $U_i$ is to valid
\begin{equation*}
\frac{\partial \hat{I}}{\partial U_i}=0, \quad i=2,\dots,n+1,
\end{equation*}
 equivalently,
\begin{equation*}
{\cal K}_{ij}U_j={\cal F}_i, \quad i=2,\dots, n+1
\end{equation*}
or, in matrix form
\begin{equation}
{\bf {\cal K}U={\cal F}}.
\end{equation}
\\
The next step is  the determination of  column vector $\cal F$ and the matrix $\cal K$ so as to confirm that they coincide with the ones given by eqs (35) and (36), respectively.\
\subsection{Determination of  $\cal F$ }
First, the column vector $\bf \tilde{F}$ is examined. Recalling  the symmetry of the discretization (see Section 2.2), it is easy to verify that the global  shape functions enjoy the symmetry
$$
 N_i(t-s)=N_m(s),\quad,
$$
where  $m-$node  is the symmetric to the $i-$node, i.e., 
$$
m=n+2-i.
$$
Then one can easily calculate the $n-$dimensional vector  $\bf \tilde{F}$ 
$$
\tilde{F}_i=\left[ \tilde{f},N_i\right]=\int_0^t \tilde{f}(s)N_i(t-s)ds=\int_0^tf_0\delta(s-s_1)N_m(s)ds=f_0\delta_{1m},
$$
for
$$
i=2,\dots,n+1\quad\text{corresponding to}\quad
m=n,\dots,1
$$
which implies
\begin{equation}
{\bf \tilde{F}}=
\begin{bmatrix}
0\\
\vdots\\
0\\
f_0
\end{bmatrix}
=
\begin{bmatrix}
0\\
\vdots\\
0\\
mv_0
\end{bmatrix}.
\end{equation}
Similarly,  the vector $\bf F$ in virtue of eqs. (28) and (37) can be determined as follows
\begin{eqnarray}
\left[f,N_i\right]=\int_0^t f(s)N_i(t-s)ds=\int_0^t f(s)N_m(s)ds=\int_{s_{m-1}}^{s_{m+1}} f(s)N_m(s)ds\nonumber\\
=\int_{s_{m-1}}^{s_m} f(s)N^{m-1}_2(s)ds+\int_{s_m}^{s_{m+1}} f(s)N^m_1(s)ds \nonumber\\
=\int_0^{\tau_m} f(s_{m-1}+s)N^{m-1}_1(s_m-s)ds+\int_0^{\tau_m} f(s_{m-1}+s)N^m_2(s_{m+1}-s)ds\nonumber
\end{eqnarray}
\begin{equation*}
\Rightarrow\left[f,N_i\right]=\left[f,N_1^{m-1}\right]_{s_{m-1}}^{s_m}+\left[f,N_2^m\right]_{s_m}^{s_{m+1}}=F_1^{m-1}+F_2^m,\quad
\begin{array}{ll}
i=2,\dots,n+1&\\
m=n,\dots,1.&
\end{array}
\end{equation*}
Thus, recalling the relation $m=n+2-i$, one concludes
\begin{eqnarray*}
F_i=\left[f,N_i\right]=F_1^{m-1}+F_2^m=F_1^{n+1-i}+F_2^{n+2-i},\\
i=2,\dots,n+1,\quad m=n,\dots,1
\end{eqnarray*}
or
\begin{equation*}
{\bf F}=
\begin{bmatrix}
F_1^{n-1}+F_2^n\\
\vdots\\
F_1^1+F_2^2\\
F_2^1
\end{bmatrix}
\end{equation*}
and by turn, $\cal F$ takes its final form
\begin{equation}
{\cal F}=
{\bf F}+\tilde{\bf F}=
\begin{bmatrix}
F_1^{n-1}+F_2^n\\
\vdots\\
F_1^1+F_2^2\\
F_2^1
\end{bmatrix}
+
\begin{bmatrix}
0\\
\vdots\\
0\\
mv_0
\end{bmatrix},
\end{equation}
thus, apart from the first row, it coincides  with eq. (36) of the previous section.
\subsection{Determinatioon of $\cal K$}
If $m$ is the symmetric to the $j$ node,  the following relation holds
\begin{equation}
[N_i,N_j]=\int_0^tN_i(s)N_j(t-s)ds=\int_0^tN_i(s)N_m(s)ds=(N_i,N_m),
\end{equation}
where $m=n+2-j$. Thus for  $j$ running from 2 to $n+1$, $m$ runs from $n$ to $1$, consequently, one can write \\
\begin{eqnarray}
{\bf K}&=&k
\begin{bmatrix}
[N_2,N_2]&[N_2,N_3]&\dots&[N_2,N_{n+1}]\\
\vdots&\vdots&\dots&\vdots\\
[N_{n+1},N_2]&[N_{n+1},N_3]&\dots&[N_{n+1},N_{n+1}]
\end{bmatrix}\nonumber\\
\nonumber\\
&=&
k\begin{bmatrix}
(N_2,N_n)&(N_2,N_{n-1})&\dots&(N_2,N_1)\\
\vdots&\vdots&\dots&\vdots\\
(N_{n+1},N_n)&(N_{n+1},N_{n-1})&\dots&(N_{n+1},N_1)
\end{bmatrix}
\end{eqnarray}
Notice that $(N_i,N_m)$ may be non--zero only if  the $i$  and  $m$ nodes coincide   or  are neighbouring, thus for a certain $i$ there are only  three cases for non-vanishing terms. Consider for instance  the case $i=m$, that is, $i$ is the symmetric to $j$ node
\begin{eqnarray*}
&&[N_i,N_j]= (N_i,N_i )\quad  (\text{for}\  i=n+2-j)\\
&&=\int_0^t{N_i^2(s)}ds=\int_{s_{i-1}}^{s_{i+1}}{N_i^2(s)}ds\\
&&=\int_{s_{i-1}}^{s_i}N_i^2(s)ds+\int_{s_i}^{s_{i+1}}N_i^2(s)ds\\
&&=\int_{s_{i-1}}^{s_i}\left(N_2^{i-1}\right)^2(s)ds+\int_{s_i}^{s_{i+1}}\left(N_1^i\right)^2(s)ds\quad  (\text{due to eq. (37)})\\
&&=\left(N_2^{i-1},N_2^{i-1}\right)+\left(N_1^i,N_1^i\right)\\
&&=\left[N_2^{i-1},N_1^{i-1}\right]_{s_{i-1}}^{s_i}+\left[N_1^i,N_2^i\right]_{s_i}^{s_{i+1}}\quad \text{ (due to eq. (28))}
\end{eqnarray*}
Thus, recalling eqs (24 ) and (42), one concludes in 
\begin{equation}
K_{ij}=K_{21}^{i-1}+K_{12}^i,\quad \text{$i$ symmetric to $j$}
\end{equation}
By the same manner, one can examine the cases $i=m-1$ and $i=m+1$,  so as to prove the general relation
\begin{equation}
K_{ij}=\left\{
\begin{array}{ll}
K_{11}^{i-1},&\quad  i=n+1-j\\
K_{21}^{i-1}+K_{12}^i,&\quad i=n=2-j\\
K_{22}^i,&\quad i=n+3-j\\
0,&\quad \text{otherwise}
\end{array}
\right.
\end{equation}
and
\begin{equation}
M_{ij}=\left\{
\begin{array}{ll}
M_{21}^{m-1}+M_{12}^m,&\quad i=m\\
M_{11}^{m-1},&\quad  i=m-1\\
M_{22}^{m},&\quad i=m+1\\
0,&\quad \text{otherwise.}
\end{array}
\right.
\end{equation}
Thus, according to eq. (41), the matrix $\cal K$ is given by the relation 
\begin{equation}
{\cal K}_{ij}=\left\{
\begin{array}{ll}
{\cal K}_{21}^{m-1}+{\cal K}_{12}^m,&\quad i=m\\
{\cal K}_{11}^{m-1},&\quad  i=m-1\\
{\cal K}_{22}^{m},&\quad i=m+1\\
0,&\quad \text{otherwise.}
\end{array}
\right.
\end{equation}
or
$$
{\cal K}_{ij}=
\begin{bmatrix}
0&\cdots&\cdots&\cdots&{\cal K}_{2(n-1)}&{\cal K}_{2n}&{\cal K}_{2(n+1)}\\
\vdots&\cdots&0&{\cal K}_{3(n-2)}&{\cal K}_{3(n-1)}&{\cal K}_{3n}&0\\
\vdots&\cdots&\cdots&\cdots&\cdots&\cdots&\vdots\\
0&\cdots&\cdots&\cdots&\cdots&\cdots&\vdots\\
{\cal K}_{(n-1)2}&{\cal K}_{(n-1)3}&{\cal K}_{(n-1)4}&0&\cdots&\cdots&\vdots\\
{\cal K}_{n2}&{\cal K}_{n3}&0&\cdots&\cdots&\cdots&\vdots\\
{\cal K}_{(n+1)2}&0&\cdots&\cdots&\cdots&\cdots&0
\end{bmatrix}.
$$
The latter, accounting for eq. (51) takes its final form
\begin{equation}
{\cal K}_{ij}=
\begin{bmatrix}
0&\cdots&\cdots&\cdots&{\cal K}^2_{11}&{\cal K}^1_{21}+{\cal K}^2_{12}&{\cal K}^1_{22}\\
\vdots&\cdots&0&{\cal K}^3_{11}&{\cal K}^2_{21}+{\cal K}^3_{12}&{\cal K}^2_{22}&0\\
\vdots&\cdots&\cdots&\cdots&\cdots&\cdots&\vdots\\
0&\cdots&\cdots&\cdots&\cdots&\cdots&\vdots\\
{\cal K}^{n-1}_{11}&{\cal K}^{n-2}_{21}+{\cal K}^{n-1}_{12}&{\cal K}^{n-2}_{22}&0&\cdots&\cdots&\vdots\\
{\cal K}^{n-1}_{21}+{\cal K}^n_{12}&{\cal K}^{n-1}_{22}&0&\cdots&\cdots&\cdots&\vdots\\
{\cal K}^n_{22}&0&\cdots&\cdots&\cdots&\cdots&0
\end{bmatrix}
\end{equation}
which is identical to the matrix given by eq, (35), if one omits the first row and the first column.
\section{A One-step Marching  Scheme}
Returnig to subesction 3.1, notice that the local linear system (26) may provide a recurrent algorithm for the computaion of the field $u$ and of its derivative $u'$. Indeed, consider the system corresponding to  the $e-$element  given by
\begin{eqnarray}
{\cal K}^{e}_{11}U_{e}+{\cal K}^{e}_{12}U_{e+1}=F^{e}_1-mV_{e+1},\nonumber\\
\\
{\cal K}^{e}_{21}U_{e}+{\cal K}^{e}_{22}U_{e+1}=F^{e}_2+mV_{e},\nonumber
\end{eqnarray}
where $V_e$ denotes the velocity at the time $s_e$, i.e., it holds
\begin{equation}
V_e=u'(s_e),\quad e=1,\dots,n.
\end{equation}
The above system can be  rewritten as
\begin{eqnarray}
V_{e+1}=\frac{1}{m}\left(-{\cal K}^{e}_{11}U_{e}-{\cal K}^{e}_{12}U_{e+1}+F^{e}_1\right),\nonumber\\
\\
U_{e+1}=\frac{1}{{\cal K}^{e}_{22}}\left({\cal K}^{e}_{21}U_{e}+mV_{e}+F^{e}_2\right).\nonumber
\end{eqnarray} 
Thus, if the displacement and velocity are known at time $s_e$, one can compute $U_{e+1}$ from eq. (55b) and, in turn, $ V_{e+1}$ from eq. (55a).  Besides, the approximation scheme can be written in the form of an implicit recurrent formula
$$
\begin{bmatrix}
{\cal K}^{e}_{12}&m\\
{\cal K}^{e}_{22}& 0
\end{bmatrix}
\begin{bmatrix}
U_{e+1}\\
V_{e+1}
\end{bmatrix}
=
\begin{bmatrix}
-{\cal K}^{e}_{11}&0\\
-{\cal K}^{e}_{21}& m
\end{bmatrix}
\begin{bmatrix}
U_{e}\\
V_{e}
\end{bmatrix}
+
\begin{bmatrix}
{ F}^{e}_1\\
{ F}^{e}_2
\end{bmatrix}
$$
If one consider a uniform discretization, i.e., fixed time step, then the local stiffness matrix does not depend on the particular element and the approximation scheme takes the form
\begin{equation}
{\bf A}{\bf W}^{I+1}={\bf B}{\bf W}^{I}+{\bf F}^{I},\quad I\in\mathbb{Z},
\end{equation}
where
\begin{equation*}
{\bf A}=\begin{bmatrix}
{\cal K}^e_{12}&m\\
{\cal K}^e_{22}& 0
\end{bmatrix},
\quad
{\bf B}=\begin{bmatrix}
-{\cal K}^e_{11}&0\\
-{\cal K}^e_{21}& m
\end{bmatrix},
\quad
{\bf  W}^I=\begin{bmatrix}
U_I\\
V_I
\end{bmatrix}, \quad
{\bf  F}^I=\begin{bmatrix}
F^I_1\\
F^I_2
\end{bmatrix}.
\end{equation*}
It is noted that the superscript $e$ does not mean dependence of ${\cal K}^e$ on the element, i.e., on the time step. On the contrary  ${\cal K}^e$ is computed once and for all elements. The column vector ${\bf F}^I$, however,  does depend on the time step and must be calculated at each time step  from the given external force. Particularly, for a free vibration problem, where the external force ${\bf F}^I$ vanishes (or is constant), the recurrent formula (56) is completely independent of the element and the evolution of the system $\bf W$ will be given by
\begin{equation}
{\bf W}^{I+1}={\bf C}^{(I+1)}{\bf W}^1,
\end{equation}
where
$$
{\bf C}={\bf A}^{-1}{\bf B}
$$
and ${\bf W}^1$ contains the initial data. To proceed further, one has to compute the matrix ${\cal K}^e$.  Recalling eqs. (20), the local shape functions and their derivatives are written
$$
N^e_1=\frac{s_{e+1}-s}{\tau}, \quad N^e_2=\frac{s-s_e}{\tau}
$$
and 
$$
{N^e_1}'=-\frac{1}{\tau}, {\quad N^e_2}'=\frac{1}{\tau},
$$
where $\tau$ is the time step. Thus, the matrices ${\bf M}^e$ and  ${\bf K}^e$ becomes 
$$
{\bf M}^e=\frac{m}{\tau}\begin{bmatrix}
1&& -1\\
-1&& 1
\end{bmatrix},
\quad
{\bf K}^e=\frac{k\tau}{3}\begin{bmatrix}
1/2&& 1\\
1&& 1/2
\end{bmatrix},
$$
respectively, which in combination (eq. (23)) result in 
\begin {equation*}
{\cal K}^e=\begin{bmatrix}
\frac{m}{\tau}+\frac{k\tau}{6}&& -\frac{m}{\tau}+\frac{k\tau}{3}\\
-\frac{m}{\tau}+\frac{k\tau}{3}&& \frac{m}{\tau}+\frac{k\tau}{6}
\end{bmatrix}
\end{equation*}
and in turn $\bf A$ and $\bf B$ take the particular form
\begin {equation}
{\bf A}=\begin{bmatrix}
-\frac{m}{\tau}+\frac{k\tau}{3}&& m\\
\frac{m}{\tau}+\frac{k\tau}{6}&& 0
\end{bmatrix},
\quad
{\bf B}=\begin{bmatrix}
-\frac{m}{\tau}-\frac{k\tau}{6}&& 0\\
\frac{m}{\tau}-\frac{k\tau}{3}&& m
\end{bmatrix}.
\end{equation}
The next question concerns the stability of the aproximation scheme (57)--(58). To secure the stability, one needs complex values of the generalized eigenvalue problem
$$
\text{det}\left({\bf B}-\lambda{\bf A}\right)=0,
$$
or
\begin{equation}
 \lambda^2-2\lambda\frac{M-K}{M+K/2}+1=0,
\end{equation}
where
\begin{equation}
M=m/\tau \quad \text{and}\quad K=k\tau/3.
\end{equation}
Requiring negative value for the determinant of eq. (65), one obtains
$$
4\left[\left(\frac{M-K}{M+K/2}\right)^2-1\right]<1\Rightarrow K/M<4,
$$
or equivalently, on account of eqs (60),
\begin{equation}
\tau^2<12m/k\quad\text{or}\quad  (\tau\omega)^2<12,
\end{equation}
where $\omega=\sqrt{k/m}$ is the frequency of the free oscillation. It can be rewritten in terms of the period $T=2\pi/\omega$ as 
\begin{equation}
\tau<\frac{\sqrt{12}}{2\pi} T\approx 0.551 T.
\end{equation}
If the time step, $\tau$ fulfils the above inequaliy, then  it is easy to confirm that for both conjugate eigenvalues, it holds 
$$
|\lambda_1|=|\lambda_2|=1.
$$
Thus, Eq. (61) (or eq. (62)) is the required condition for the numerical stability of the approximation scheme.
\section{Computational Applications}
In this section, some computational examples are presented to support the theoretical and numerical  considerations of the previous sections. 
Uniform discretizations (elements of equal length) are used so that the results of the two approximation schemes can be fully comparable.  In that case the local stifness matrices are fixed  for every element (see Appendix A)
$$
{\bf K}^e=\frac{m{\tau}}{3}\begin{bmatrix}
1/2&& 1\\
1&& 1/2
\end{bmatrix},
\quad
{\bf M}^e=\frac{m}{{\tau}}\begin{bmatrix}
1&& -1\\
-1&& 1
\end{bmatrix}.
$$
Taking into account  eqs. (74) and (75) of Appendix B, the  global matrices for a uniform discretization become 
\begin{equation*}
{\bf K}=
\frac{k\tau}{3}\begin{bmatrix}
0 & \dots &\dots & 1/2 &2&1/2\\
0 & \dots & 1/2&2&1/2& 0\\
\vdots &\vdots   & \vdots &\vdots & \vdots & \vdots\\
 2 &1/2 &0 &0&\dots&0\\
 1/2 & 0 & 0 &\dots&\dots&0
\end{bmatrix},
\quad
{\bf M}=
\frac{m}{\tau}\begin{bmatrix}
0 & \dots &\dots &-1 &2&-1\\
0 & \dots &-1&2&-1& 0\\
\vdots &\vdots   & \vdots &\vdots & \vdots & \vdots\\
 2&-1&0 &0&\dots&0\\
-1 & 0 & 0 &\dots&\dots&0
\end{bmatrix}.
\end{equation*}
All examples  run with the constants  $k=9$ and $m=1$, so as the natural frequency of the system to become $\omega=3$. \\
\\
First, both  approximation schemes are tested  for the free vibration problem. In that case, the external force vanishes, i.e.,  $f=0$ and the system is excited by the initial velocity or the initial position. For the examples presented here, we have taken $u(0)=0$ and $u'(0)=2$. In Fig. 1, the values of the position versus time  are shown for three distinct values of the element length (or time step). Simultaneously, comparison with the analytical solution  is demonstrated in each case. Analogous results and comparisons  are presented in Table 1.\\
\\
The second example concerns a forced vibration problem. For this reason, an excitation term of the form 
$$
f(s)=f_0\ \text{sin}(\Omega t)
$$
was taken, where  $\Omega$ is the frequency and $f_0$  the amplitude of the external force.  The nodal forces, i.e. eqs (25a) can be analytically calculated as  (see Appendix C)
\begin{eqnarray}
F_1^e&=&-(f_0/\Omega) \cos (\Omega s_{e+1})+\left(f_0/(\tau_e\Omega^2)\right)\left[\sin(\Omega s_{e+1})-\sin(\Omega s_e)\right] \nonumber\\
F_2^e&=&(f_0/\Omega) \cos (\Omega s_{e})-\left(f_0/(\tau_e\Omega^2)\right)\left[\sin(\Omega s_{e+1})-\sin(\Omega s_e)\right]\nonumber
\end{eqnarray}
In Fig. 2 and Table 2, the responce of the system to the external excitation is demonstraded for various values of $\tau$ and is compared with the analytical solution. Also, in Fig. 3a, the algorithms for the two  approximation solutions were run for a longer time to reveal the characteristic pattern of the forced vibration and   the particular case where the exciation frequency meets the natural velocity of the system is shown in Fig. 3b.\\
\\
It is remarkable that both approximation schemes, the FEM and the One--step produce the same numbers, at least at the accuracy of four decimal digits which are used here. To find differences between the two schemes, one has to proceed further accounting for at least twelve additional decimal digits, as one can confirm by looking at the Table 3. \\  
\section{Conclusions}
In \cite{KALPAKIDES-et-al2019}, a generalized Hamilton's principle has been developed which is  free of the main  inconsistency  of the classical Hamilton's principle, i.e., its weakness to account properly for  the initial conditions of a dynamical problem. The key role to this developement is the use of the convolution instead of the $L^2$ inner product.\\
\\
In the present paper, this generalized principle is exploited to formulate new approximation schemes for a simple initial value problem.  More specifically, a consistent finite element formulation is fully developed. The use of convolution introduces significant changes in the finite dimensional version of the problem.  The  stiffness matrix is now  symmetric with respect to its second diagonal and the load vector is inversed.  Besides, a recurrent formula stems from the consideration of the algebraic system at the element which provides an one--step approximation scheme, the numerical  stability of which is examined. Both approximation schemes are tested by computational examples concerning the free and forced vibration of an  one degree of freedom oscillator.\\
\\
It seems that the analysis presented here may be applied to initial value problems in general. Actually, the essential conclusion of this work is the  possibility to formulate FE schemes for initial value problems, provided that an appropriate convolutional variational formulation has been previously developed.\\
\\
The next challenge is to further proceed  to the formulation of  FE scheme for a unified space--time continuoum, i.e., to an undivided  FE scheme for  initial--boundary value problems.
\newpage
\appendix
{\LARGE{\bf APPENDICES}}
\section{ Calculation of Local Matrices and Forces}
The terms of eq. (21) are calculated one by one:
\begin{eqnarray}
\bullet&& [{u^e}',{u^e}']^{s_{e+1}}_{s_e}= \int_0^{\tau_e}{u^e}'(s_e+s){u^e}'(s_{e+1}-s)ds=\nonumber\\
&& {u^e_1}^2 \int_0^{\tau_e}{N^e_1}'(s_e+s){N^e_1}'(s_{e+1}-s)ds+u^e_1u^e_2\int_0^{\tau_e}{N^e_1}'(s_e+s){N^e_2}'(s_{e+1}-s)ds\nonumber \\
&&+u^e_2u^e_1\int_0^{\tau_e}{N^e_2}'(s_e+s){N^e_1}'(s_{e+1}-s)ds+ {u^e_2}^2 \int_0^{\tau_e}{N^e_2}'(s_e+s){N^e_2}'(s_{e+1}-s)ds \nonumber\\
&&={u^e_1}^2 \left[{N^e_1}',{N^e_1}'\right]^{s_{e+1}}_{s_e}+2u^e_1u^e_2\left[{N^e_1}',{N^e_2}'\right]^{s_{e+1}}_{s_e}
+{u^e_2}^2 \left[{N^e_2}',{N^e_2}'\right]^{s_{e+1}}_{s_e}\nonumber
\end{eqnarray}
\begin{equation}
\Rightarrow\frac{m}{2}[{u^e}',{u^e}']^{s_{e+1}}_{s_e}=\frac{1}{2}\left({u^e_1}^2M^e_{11}+2 u^e_1 u^e_2M^e_{12}+{u^e_2}^2M^e_{22}\right)=\frac{1}{2}{{\bf u}^e}^T{\bf M}^e{\bf u}^e,
\end{equation}
where
\begin{equation}
{\bf M}^e=M^e_{ij}=m\left[{N^e_i}',{N^e_j}'\right]^{s_{e+1}}_{s_e},\quad
\begin{array}{ll}
i=1,2,&\\
e=1,\dots,n.&
\end{array}
\end{equation}
and 
\begin{equation}
{\bf M}^e=\frac{m}{{\tau}^e}\begin{bmatrix}
1&& -1\\
-1&& 1
\end{bmatrix}.
\end{equation}
\begin{eqnarray}
\bullet&& [u^e,u^e]^{s_{e+1}}_{s_e}=\int_0^{\tau_e}u^e(s_e+s)u^e(s_{e+1}-s)ds=\nonumber\\
&& {u^e_1}^2 \int_0^{\tau_e}N^e_1(s_e+s)N^e_1(s_{e+1}-s)ds+u^e_1u^e_2\int_0^{\tau_e}N^e_1(s_e+s)N^e_2(s_{e+1}-s)ds\nonumber \\
&&+u^e_2u^e_1\int_0^{\tau_e}N^e_2(s_e+s)N^e_1(s_{e+1}-s)ds +{u^e_2}^2 \int_0^{\tau_e}N^e_1(s_e+s)N^e_1(s_{e+1}-s)ds \nonumber\\
&&={u^e_1}^2 \left[N^e_1,N^e_1\right]^{s_{e+1}}_{s_e}+2u^e_1u^e_2\left[N^e_1,{N^e_2}'\right]^{s_{e+1}}_{s_e}
+{u^e_2}^2 \left[N^e_2,N^e_2\right]^{s_{e+1}}_{s_e}\nonumber
\end{eqnarray}
\begin{equation}
\Rightarrow\frac{k}{2} [u^e,u^e]^{s_{e+1}}_{s_e}=\frac{1}{2}\left({u^e_1}^2K^e_{11}+2 u^e_1 u^e_2K^e_{12}+{u^e_2}^2K^e_{22}\right)
=\frac{1}{2}{{\bf u}^e}^T{\bf K}^e{\bf u}^e,
\end{equation}
where
\begin{equation}
{\bf K}^e=K^e_{ij}= k\left[N^e_i,N^e_j\right]^{s_{e+1}}_{s_e},\quad
\begin{array}{ll}
 i,j=1,2,&\\
\ e=1,\dots, n&
\end{array}
\end{equation}
and 
\begin{equation}
{\bf K}^e=\frac{m{\tau}^e}{3}\begin{bmatrix}
1/2&& 1\\
1&& 1/2
\end{bmatrix}.
\end{equation}
\begin{eqnarray}
\bullet &&[f,u^e]^{s_{e+1}}_{s_e}=\int_0^{\tau_e}f(s_e+s)u^e(s_{e+1}-s)ds=\nonumber\\
&&\int_0^{\tau_e}f(s_e+s)\left(u^e_1N^e_1(s_{e+1}-s)+u^e_2N^e_2(s_{e+1}-s)\right)ds=\nonumber\\
&&u^e_1\left(\int_0^{\tau_e}f(s_e+s)N^e_1(s_{e+1}-s)ds\right)+u^e_2\left(\int_0^{\tau_e}f(s_e+s)N^e_2(s_{e+1}-s)ds\right)\nonumber
\end{eqnarray}
\begin{equation}\Rightarrow[f,u^e]^{s_{e+1}}_{s_e}= u^e_1F^e_1+u^e_2 F^e_2={{\bf u}^e}^T{\bf F}^e,
\end{equation}
where
$$
{\bf F}^e=F^e_i=\int_0^{\tau_e}f(s_e+s)N^e_i(s_{e+1}-s)ds,\ i=1,2,\ e=1,\dots,n
$$
\begin{eqnarray}
\bullet\qquad  \left[\tilde{f^e_2},u^e\right]^{s_{e+1}}_{s_e}=\int_0^{\tau_e}\tilde{f}^e_2(s_e+s)u^e(s_{e+1}-s)ds=
\int_0^{\tau_e}f^e_2\delta(s-s_{e+1})u^e(s_{e+1}-s)ds\nonumber
\end{eqnarray}
\begin{equation}
\Rightarrow\left[\tilde{f^e_2},u^e\right]^{s_{e+1}}_{s_e}=f^e_2u^e_1
\end{equation}
\begin{eqnarray}
\bullet\qquad  \left[\tilde{f^e_1},u^e\right]^{s_{e+1}}_{s_e}=\int_0^{\tau_e}\tilde{f}^e_1(s_e+s)u^e(s_{e+1}-s)ds=
\int_0^{\tau_e}f^e_1\delta(s-s_e)u^e(s_{e+1}-s)ds\nonumber
\end{eqnarray}
\begin{equation}
\Rightarrow\left[\tilde{f^e_1},u^e\right]^{s_{e+1}}_{s_e}=f^e_1u^e_2
\end{equation}
\section{ Calculation of $[N_i,N_j]$ and $[{N_i}',{N_j}']$ }
The $n+2-j$  node is symmetric to the $j$ node. The integral
$$
[N_i,N_j]=(N_i,N_{n+2-j})=\int_0^tN_i(s)N_{n+2-j}(s)ds
$$
may be non-zero only if the $i$  and  $n+2-j$ nodes coincide  or are neighbouring, that is, for a particular $i$  there are only  three cases for non-vanishing terms\\
$\bullet$ {\em Case 1.}\quad  $n+2-j= i$
\begin{eqnarray*}
[N_i,N_j]&=&\int_0^tN_i(s)N_j(t-s)ds=\int_0^tN_i(s)N_{n+2-j}(s)ds=\int_0^tN_i(s)N_i(s)ds\\
&=&\int_{s_{i-1}}^{s_{i+1}}{N_i(s)}^2ds =\int_{s_{i-1}}^{s_i}\left(\frac{s-s_{i-1}}{\tau_{i-1}}\right)^2ds+\int_{s_{i}}^{s_{i+1}}\left(\frac{s_{i+1}-s}{\tau_{i}}\right)^2ds\\
&=&\frac{\tau_{i-1}}{3}+\frac{\tau_{i}}{3}\quad\Rightarrow [N_i,N_j]=\frac{\tau_i+\tau_{i+1}}{3}.
\end{eqnarray*}
$\bullet$ {\em Case 2.}\quad $n+2-j= i+1$
\begin{eqnarray*}
[N_i,N_j]=\int_0^tN_i(s)N_{n+2-j}(s)ds=\int_0^tN_i(s)N_{i+1}(s)ds=\int_{s_i}^{s_{i+1}}N_i(s)N_{i+1}(s)ds=\\
\int_{s_i}^{s_{i+1}}\left(\frac{s_{i+1}-s}{\tau_i}\right)\left(\frac{s-s_i}{\tau_i}\right)ds=\frac{\tau_i}{6}\Rightarrow [N_i,N_j]=\frac{\tau_i}{6}.
\end{eqnarray*}
$\bullet$ {\em Case 3.}\quad $n+2-j= i-1$
\begin{eqnarray*}
[N_i,N_j]=\int_0^tN_i(s)N_{n+2-j}(s)ds=\int_0^tN_i(s)N_{i-1}(s)ds=\int_{s_{i-1}}^{s_{i}}N_i(s)N_{i-1}(s)ds=\\
\int_{s_{i-1}}^{s_{i}}\left(\frac{s-s_{i-1}}{\tau_{i-1}}\right)\left(\frac{\tau_i-s}{\tau_{i-1}}\right)ds=\frac{\tau_{i-1}}{6}\Rightarrow [N_i,N_j]=\frac{\tau_{i-1}}{6}.
\end{eqnarray*}
Thus, one concludes in the general formula
\begin{equation}
[N_i,N_j]=\left\{
\begin{array}{cl}
\tau_i/6,&\quad i=n+1-j\\
(\tau_{i-1}+\tau_i)/3,&\quad  i=n+2-j\\
\tau_i/6,&\quad i=n+3-j\\
0,&\quad \text{otherwise}
\end{array}.
\right.
\end{equation}
Analogously, one may prove
\begin{equation}
[N'_i,N'_j]=\left\{
\begin{array}{cl}
-1/\tau_i,&\quad i=n+1-j\\
1/\tau_{i-1}+1/\tau_i,&\quad  i=n+2-j\\
1/\tau_{i-1},&\quad i=n+3-j\\
0,&\quad \text{otherwise.}
\end{array}.
\right.
\end{equation}
After that, the matrices ${\bf K}$  and ${\bf M}$ are written
\begin{equation}
{\bf K}=
\frac{k}{3}\begin{bmatrix}
0 & \dots &\dots & \tau_2/2 &\tau_1+\tau_2&\tau_1/2\\
0 & \dots & \tau_3/2&\tau_2+\tau_3&\tau_2/2& 0\\
\vdots &\vdots   & \vdots &\vdots & \vdots & \vdots\\
 \tau_{n-1}+\tau_n &\tau_{n-1}/2 &0 &0&\dots&0\\
 \tau_n/2 & 0 & 0 &\dots&\dots&0
\end{bmatrix}
\end{equation}
and
\begin{equation}
{\bf M}=
m\begin{bmatrix}
0 & \dots &\dots &-1/ \tau_2 &1/\tau_1+1/\tau_2&-1/\tau_1\\
0 & \dots &-1/ \tau_3&1/\tau_2+1/\tau_3&-1/\tau_2& 0\\
\vdots &\vdots   & \vdots &\vdots & \vdots & \vdots\\
 1/\tau_{n-1}+1/\tau_n &-1/\tau_{n-1} &0 &0&\dots&0\\
-1 \tau_n & 0 & 0 &\dots&\dots&0
\end{bmatrix},
\end{equation}
respectively.
\section{Calculation of Nodal Forces}
From eqs (20),  (25a) and (63), one obtains
\begin{eqnarray}
F_1^e&=&\int_0^{\tau_e}f(s_e+s)N_1^e(s_{e+1}-s)ds=f_0\int_0^{\tau_e}\cos\left(\Omega(s_e+s)\right) \frac{s_{e+1}-(s_{e+1}-s)}{\tau_e}ds\nonumber\\
&=&(f_0/\tau_e)\int_0^{\tau_e}\cos\left(\Omega(s_e+s)\right)sds.
\end{eqnarray}
Similarly
\begin{eqnarray}
F_2^e&=&\int_0^{\tau_e}f(s_e+s)N_2^e(s_{e+1}-s)ds=f_0\int_0^{\tau_e}\cos\left(\Omega(s_e+s)\right)\frac{(s_{e+1}-s)-s_e}{\tau_e}ds\nonumber\\
&=&(f_0/\tau_e)\int_0^{\tau_e}\cos\left(\Omega(s_e+s)\right)(\tau_e-s)ds.
\end{eqnarray}
Also, it is easy to see that 
$$
\int_0^{\tau_e}\cos\left(\Omega(s_e+s)\right)sds=(1/\Omega)\left[\sin(\Omega(s_e+\tau_e))\tau_e+(1/\Omega)\left( \cos (\Omega s_{e+1})-\cos (\Omega s_e)\right)\right],
$$
with the aid of which eqs (76) and (77) become
\begin{eqnarray}
F_1^e=(f_0/\Omega)\sin(\Omega s_{e+1})+(f_0/\tau_e \Omega^2)\left[\cos(\Omega s_{e+1})-\cos(\Omega s_e)\right],\\
F_2^e=-(f_0/\Omega)\sin(\Omega s_e)-(f_0/\tau_e \Omega^2)\left[\cos(\Omega s_{e+1})-\cos(\Omega s_e)\right].
\end{eqnarray}
\newpage
\bibliographystyle{plain}
\bibliography{FEM7}
\newpage
\begin{figure}[p]
\centering
\includegraphics[scale=0.45]{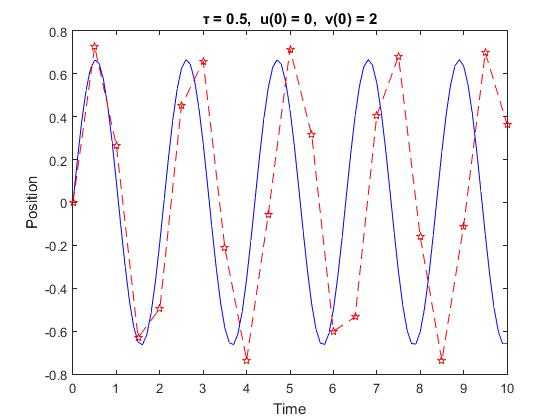}
\includegraphics[scale=0.45]{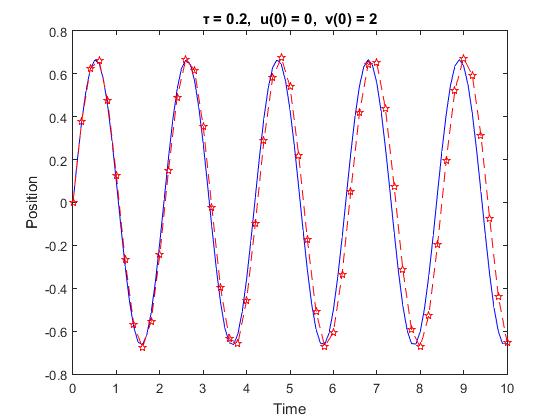}
\includegraphics[scale=0.45]{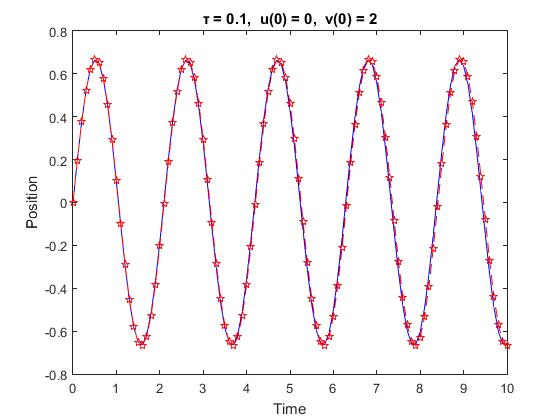}
\caption{FEM approximation for the free vibration for different values of the $\tau$, the continuous line is the exact solution}
\end{figure}

\begin{figure}[p]
\centering
\includegraphics[scale=0.45]{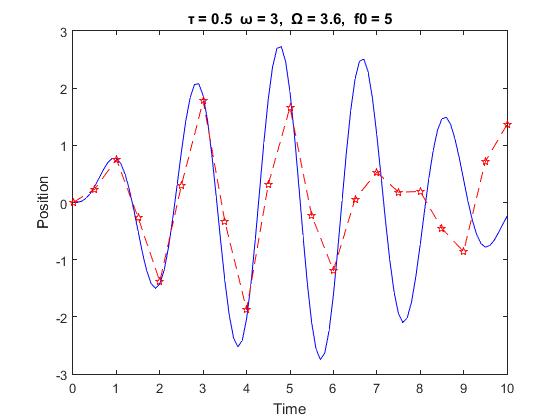}
\includegraphics[scale=0.45]{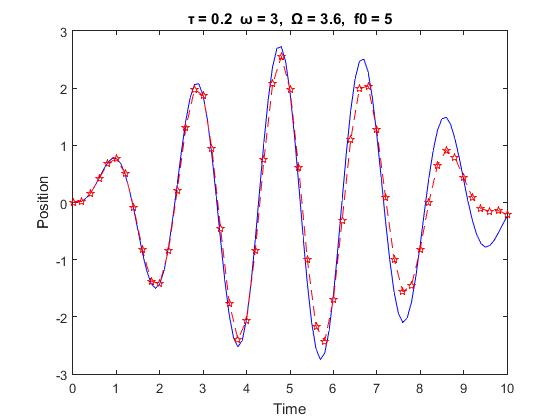}
\includegraphics[scale=0.45]{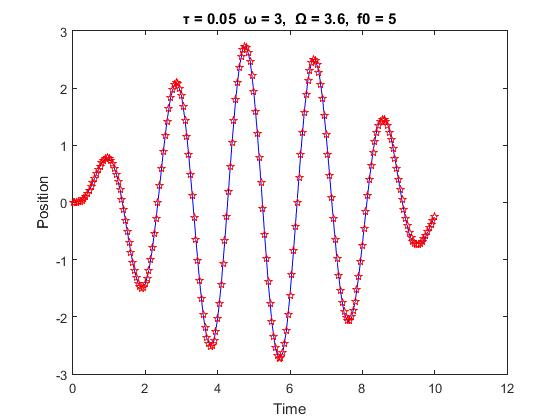}
\caption{FEM and one--step approximation schemes (asterisks) for the forced vibration, for different values of the $\tau$ versus  the exact solution (blue  solid line)}
\end{figure}
\begin{figure}[p]
\includegraphics[scale=0.45]{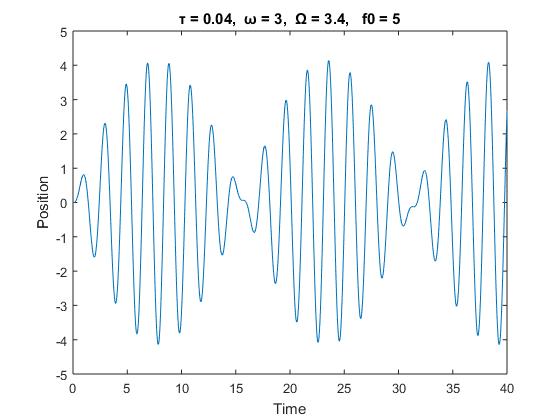}
\includegraphics[scale=0.45]{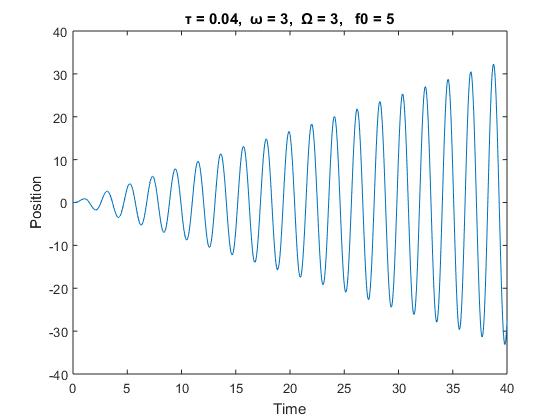}
\caption{On the left: Characteristic pattern of a vibration imposed by external force. On the right: Resonance}
\end{figure}
\begin{table}
\begin{tabular}{|c| r r r r r r| r|}
\hline
  &  &  &FEM, & One--step & &  &   Exact\\
\hline
Time  &$ \tau=0.1$&$ \tau=0.05$ & $\tau= 0.025$& $\tau=  0.02$&  $\tau=0.0125$  & $ \tau= 0.01$  &   \\
\hline
    1&    0.1018&    0.0960&    0.0946  &  0.0944   & 0.0942 &   0.0942&    0.0941\\
    2 &  -0.2012  & -0.1900  & -0.1872 &  -0.1869&   -0.1865 &  -0.1864 &  -0.1863\\
    3  &  0.2960 &   0.2801 &   0.2761  &  0.2756 &   0.2751 &   0.2750  &  0.2747\\
    4  & -0.3839 &  -0.3643 &  -0.3594  & -0.3588 &  -0.3581&   -0.3580 &  -0.3577\\
    5   & 0.4628  &  0.4410  &  0.4354    &0.4347  &  0.4340  &  0.4338   & 0.4335\\
    6  & -0.5309  & -0.5085  &-0.5026  & -0.5019  & -0.5012  & -0.5010  & -0.5007\\
    7  &  0.5867  &  0.5654    &0.5597   & 0.5590  &  0.5583  &  0.5581  &  0.5578\\
    8  & -0.6288  & -0.6105  & -0.6054  & -0.6048  & -0.6042  & -0.6040  & -0.6037\\
    9   & 0.6563  &  0.6429 &    0.6390  &  0.6385  &  0.6379  &  0.6378  &  0.6376\\
   10   &-0.6685  & -0.6619 &  -0.6595   &-0.6592 &  -0.6589  & -0.6588  & -0.6587\\
\hline
\end{tabular}
\caption{Approximate solution of the free vibration problem for various time steps.}
\end{table}
\begin{table}
\begin{tabular}{|c| r r r r r r| r|}
\hline
  &  &  &FEM, & One--step & &  &   Exact\\
\hline
Time  &$ \tau=0.1$&$ \tau=0.05$ & $\tau= 0.025$& $\tau=  0.02$&  $\tau=0.0125$  & $ \tau= 0.01$  &   \\
\hline
    1&  0.7713 &   0.7722&    0.7725&    0.7725 &   0.7725&    0.7725&    0.7726\\
    2&   -1.4247&   -1.4253&   -1.4254&   -1.4254&   -1.4255&   -1.4255&   -1.4255\\
    3&    1.8657&    1.8638&    1.8632&    1.8631&    1.8630 &   1.8630&    1.8630\\
    4&   -2.0418&   -2.0349&   -2.0329&   -2.0327&   -2.0324&   -2.0324&   -2.0323\\
    5&    1.9523&    1.9386&    1.9348&    1.9343&    1.9338&    1.9337&    1.9335\\
    6&   -1.6478&   -1.6270&   -1.6212&   -1.6205&   -1.6197&   -1.6196&   -1.6192\\
    7&    1.2188&    1.1927&    1.1853&    1.1844&    1.1834&    1.1832&    1.1828\\
    8&   -0.7766&   -0.7493&   -0.7413&   -0.7403&   -0.7392&   -0.7390&   -0.7385\\
    9&    0.4296&    0.4072&    0.4001&    0.3992&    0.3982&    0.3980 &   0.3976\\
   10&   -0.2607&   -0.2506&   -0.2464&   -0.2458&  -0.2452&   -0.2450&   -0.2448\\
\hline
\end{tabular}
\caption{Approximate values of position  for the forced vibration problem, $f_0=5,\ \Omega=3.6.$}
\end{table}
\begin{table}
\centering
\begin{tabular}{|c|| c  c c|}
\hline
   &  &$ \tau=0.5$ &   \\
\hline
time& $F$ &   $ O$& $D\times10^{14}$ \\
\hline
    1&       0.7678&    0.7678 &        0\\
    2 &     -1.4215&   -1.4215&    0.0444\\
    3&       1.8691&    1.8691&         0\\
    4&      -2.0581&   -2.0581&    0.2665\\
    5&       1.9850&    1.9850&   -0.1110\\
    6&      -1.6952&   -1.6952&    0.2220\\
    7 &      1.2718&    1.2718&   -0.0888\\
    8 &     -0.8185&   -0.8185&    0.2554\\
    9&       0.4369&    0.4369&   -0.4441\\
   10 &     -0.2060&   -0.2060&    0.5718\\
\hline
\end{tabular}
\begin{tabular}{| c c c|}
\hline
 & $\tau=0.01$& \\
\hline
  F& O&  $D\times10^{12}$  \\
\hline
      0.7725 &   0.7725&   -0.0356\\
      -1.4255&   -1.4255 &   0.0429\\
       1.8630&    1.8630 &  -0.0551\\
      -2.0324 &  -2.0324 &   0.0226\\
      1.9337 &   1.9337  & -0.0628\\
      -1.6196&   -1.6196 &   0.2720\\
       1.1832 &   1.1832 &  -0.3257\\
      -0.7390 &  -0.7390 &   0.3303\\
        0.3980&    0.3980&   -0.3147\\
     -0.2450 &  -0.2450  &  0.1495\\
\hline
\end{tabular}
\caption{Forced oscillation for  $\Omega=3.6,\ f_0=5,$  $F$: FEM approximation, $O$: One--step approximation, $D=F-O.$}
\end{table}
\end{document}